
\magnification=\magstep1
\input amstex
\documentstyle{amsppt}
\leftheadtext{K. J. B\"or\"oczky, E. Makai, Jr.}
\rightheadtext{Planar Blaschke-Santal\'o inequality}
\topmatter
\title Remarks on planar Blaschke-Santal\'o inequality\endtitle
\author K. J. B\"or\"oczky, E. Makai, Jr.* \endauthor
\address Alfr\'ed R\'enyi Mathematical Institute, 
\newline
Hungarian Academy of Sciences,
\newline
H-1364 Budapest, Pf. 127, 
\newline
HUNGARY
\newline
{\rm{http://www.renyi.hu/\~{}makai}}
\endaddress
\email carlos\@renyi.hu, makai.endre\@renyi.mta.hu\endemail
\keywords Blaschke-Santal\'o inequality, stability, polygons, $n$-fold
rotational symmetry\endkeywords
\subjclass {\it 2000 Mathematics Subject Classification.} Primary:
52A40. Secondary: 52A38, 52A10\endsubjclass
\abstract We prove the Blaschke-Santal\'o inequality restricted to $n$-gons:
the extremal polygons are the affine regular $n$-gons. If either the John or
the L\"owner 
ellipse of a planar $o$-symmetric convex body $K$ is the unit circle
about $o$, then a sharpening of the Blaschke-Santal\'o inequality holds: even
the aritmetic mean $\left( V(K) + V( K^*) \right) /2 $ is at least $\pi $. 
We give stability variants of the Blaschke-Santal\'o inequality for the
plane. If for some $n \ge 3$ the planar convex body $K$ is 
$n$-fold rotationally symmetric about $o$, then we give the exact maximum of
$V(K^*)$, as a function of $V(K)$ and the area of either the John or the
L\"owner ellipse.
\endabstract
\endtopmatter\document


We introduce some notation that will be used
in the paper.
For $x, y\in R^2$, $x \ne y$ we 
write aff$\,\{ x, y \} $ for
the line passing through $x$ and $y$. 
For compact sets or points $X_1,\ldots,X_k$ in ${\Bbb R}^2$, let
$[X_1,\ldots,X_k]$ denote the convex hull of their union, where the $X_i$'s
which are points are replaced by $\{ X_i \} $.
For $x, y\in R^2$ we write $|xy|$ for
the length of the segment $[x,y]$.
For any point $p\neq o$, let $p^*$ be the polar line
with equation $\langle x, p \rangle =1$. For
$l=p^*$, let $p=l^*$. By the intersection point
of two parallel lines we mean their common point at infinity (in the
projective plane). For two non-negative 
quantities $f,g$ we write $g=\Theta (f)$ ($g$ is of exact order $f$), 
if $g=O(f)$ and $f=O(g)$.

\head Optimality of regular $n$-gons in the $o$-symmetric case\endhead

First we treat $o$-symmetric polygons because the argument is 
much easier in this case.

\proclaim{Theorem A} For even $n\geq 6$, if $K$ is
an $o$-symmetric convex polygon of at most $n$ vertices, then
$$
V(K) V(K^*)\leq n^2\sin^2\frac{\pi}n,
$$
with equality if and only if $K$ is a regular $n$-gon.
\endproclaim

\demo{Proof} {\bf{1.}}
Let  $K$ maximize
the area product among $o$-symmetric convex polygons of at
most $n$ vertices. Since $n^2 \sin ^2 ( \pi /n )$ is strictly monotonically
increasing with $n$, by induction we may suppose that $K$ has exactly $n$
sides.
By compactness of 

\newpage

affine equivalence classes of
convex bodies such a
$K$ exists. In this case, $K^*$  is also a maximal body.
The vertices of $K$ will be denoted by $x_1, x_2,\ldots,x_n$, in positive
sense.

First we observe that $K$ has exactly $n$ vertices. In fact, if $k<n$, 
consider the
vertex $x_2$ of $K$. Let $l$ be a line such that $K \cap l = \{ x_2 \} $. Let
$\varepsilon >0$ be sufficiently small. Push $l$ parallelly towards $K$ throgh
a distance $\varepsilon $. Then
the vertex $x_2$ is cut off. (Without further mentioning
we will make always the same changes at the opposite vertex
$-x_2$ so that $K$ remains $o$-symmetric.) We write $K_{{\text{new}}}$ for the
so obtained new $o$-symmetric polygon, that still has at most $n$ vertices by
hypothesis.
Then $l$ 
cuts off an area $\Theta (\varepsilon ^2)$ from $K$, so 
$$
V(K_{{\text{new}}})=
V(K) - \Theta (\varepsilon ^2)\,.
$$
Then $K^*$ has neighbouring vertices 
${\text{aff}}\, \{ x_1,x_2 \} ^*$ 
and ${\text{aff}}\,\{ x_2,x_3 \} ^*$
but $K_{{\text{new}}}^*$
will have a new vertex $l^*$ between these two old vertices.  
The distance of $l^*$ and the side $x_2^*$ of $K$ is $\Theta (
\varepsilon )$. Hence $K_{{\text{new}}}^*$ 
is obtained from $K^*$ by adding to it
a small triangle of height, and also of area $\Theta ( \varepsilon
)$. 
Thus 
$$
V(K_{{\text{new}}}^*)= V(K)+\Theta (\varepsilon ) \,.
$$
Therefore for the volume products we have
$$
V(K_{{\text{new}}})V(K_{{\text{new}}}^*)= V(K)V(K^*) +\Theta (\varepsilon )
>V(K)V(K^*)\,,
$$ 
a contradiction.

Next we prove two basic properties of $K$:

\medskip

\noindent(i) If $x_1,x_2,x_3$ are consecutive vertices of $K$,
then $o$, $(x_1+x_3)/2$ and $x_2$ are collinear. 

\medskip

\noindent(ii) If $y_1,y_2,y_3,y_4$ are consecutive vertices of $K$,
and $m$ is the intersection point of aff$\,\{ y_1, y_2 \} $ 
and aff$\,\{ y_3, y_4 \} $ (observe that $m$
is a finite point separated from $K$ by aff\,$\{ y_2,y_3 \} $),
then $o$, $(y_2+y_3)/2$ and $m$ are collinear. 

\medskip

Observe that both (i) and (ii) are affine invariant. Then we may suppose that
in (i) $x_1,x_2,x_3 \in S^1$, and in (ii) that aff\,$\{ y_1, y_2\} $,
aff\,$\{ y_2, y_3 \} $, aff\,$\{ y_3, y_4 \} $ are
tangents to $S^1$. Then (i) expresses that $\{ x_1,x_2,x_3 \} $ is symmetrical
to the axis $\Bbb{R} \cdot (x_1+x_3)/2$, while (ii) expresses that 
$\{ $aff\,$\{ y_1, y_2\} $,
aff\,$\{ y_2, y_3 \} $, aff\,$\{ y_3, y_4 \} \}$ is symmetrical to $\Bbb{R}
\cdot m$.
Otherwise said, in the metric of the $o$-symmetric ellipse containing
$x_1,x_2,x_3$, the rotation carrying $x_1$ to $x_2$ carries $x_2$ to $x_3$,
and in the metric of the $o$-symmetric ellipse touched by 
$ $aff\,$\{ y_1, y_2\} $,
aff\,$\{ y_2, y_3 \} $, aff\,$\{ y_3, y_4 \} $, the rotation carrying
$ $aff\,$\{ y_1, y_2\} $ to aff\,$\{ y_2, y_3 \} $ carries 
aff\,$\{ y_2, y_3 \} $ to aff\,$\{ y_3, y_4 \} $.
Then these two symmetry properties, one taken for $K$, the other
taken for $K^*$, clearly imply each other.

Since (i) for $K^*$ is equivalent to (ii) for $K$, therefore
it is sufficient to prove (i).

We suppose that (i) does not hold, and seek a contradiction
with a method going back to Mahler. (Recall that we assumed
$x_1,x_2,x_3 \in S^1$.)
We may assume that $x_1$ and $x_2$ lie on the
same open side of aff$\,\{ o,(x_1+x_3)/2 \} $. For $i=1,5$, $2.5$, we write
$y_i$ 

\newpage

for the intersection of $x_i^*$ and $x_2^*$,
and hence $y_{1.5}$ and $y_{2.5}$ are consecutive vertices of $K^*$. 
Let $l$ be the line
through $x_2$ parallel to aff$\,\{ x_1, x_3 \} $. 
We move $x_2$
into a position $\tilde{x}_2$ along $l$ towards  
aff$\,\{ o, (x_1+x_3)/2 \} $. 
In other words,
$\tilde{x}_2-x_2=\varepsilon (x_3-x_1)$ for small $\varepsilon>0$,
where $\varepsilon$ is small enough to ensure that
$x_1$ and $\tilde{x}_2$ lie on the
same open side of both aff$\,\{ o,(x_1+x_3)/2 \} $ 
and the line containing the other side of $K$ with endpoint $x_3$ (say,
aff\,$\{ x_3,x_4 \} $). (That is, we move toward the symmetric situation.)
Therefore there exists 
an $o$-symmetric convex polygon $\widetilde{K}$ obtained from $K$ by removing
$x_2$ and $-x_2$ from the set of vertices, 
and adding $\tilde{x}_2$ and $-\tilde{x}_2$.
Then $V(\widetilde{K})=V(K)$.

For $i=1.5$, or $i=2.5$, let
$\tilde{y}_i$ be the intersection of $x_1^*$ and $\tilde{x}_2^*$, or
$x_3^*$ and $\tilde{x}_2^*$, respectively,
and hence $\tilde{y}_{1.5}$ and $\tilde{y}_{2.5}$ are
the two new vertices of $\widetilde{K}^*$ 
replacing the vertices $y_{1.5}$ and $y_{2.5}$ of $K^*$.
In addition $\tilde{y}_{1.5} \not\in K^*$ and $\tilde{y}_{2.5}\in K^*$, and
$l^*=[y_{1.5},y_{2.5}] \cap [\tilde{y}_{1.5},\tilde{y}_{2.5}]$, moreover 
$$
\|l^*-y_{1.5}\|>\|l^*-y_{2.5}\| \,\,
{\text{ and hence }}\,\,
\|l^*-\tilde{y}_{1.5}\|>\|l^*-\tilde{y}_{2.5}\|,
$$
for $\varepsilon >0$ sufficiently small.
Since the triangles $[l^*,y_{1.5},\tilde{y}_{1.5}]$
and $[l^*,y_{2.5},\tilde{y}_{2.5}]$ have the same angle at $l^*$, we have
$$
V(\tilde{K}^*)-V(K^*)=
2 \left( V([l^*,y_{1.5},\tilde{y}_{1.5}])-V([l^*,y_{2.5},\tilde{y}_{2.5}]) 
\right) > 0
$$
(observe that because of $o$-symmetry, an  
analogous change has to be made at the
side opposite to $y_{1.5}y_{2.5}$, that gives the factor $2$). 
This contradiction
proves (i).

{\bf{2.}}
Now we prove Theorem A based on (i) and (ii).
Applying a linear transformation, we may assume that $x_1,x_2,x_3 \in S^1$,
and hence $x_2$ lies on the
perpendicular bisector of $[x_1,x_3]$ by (i). 

Now (ii) yields that aff$\,\{ x_1, x_2 \} $
and aff$\,\{ x_3, x_4 \} $ are symmetric with respect to
the perpendicular bisector of $[x_2,x_3]$ (that contains $o$). It follows
that aff$\, \{ x_2, x_3 \} $
and aff$\,\{ x_3, x_4 \} $ are symmetric with respect to
aff$\, \{ o, x_3 \} $. Together with (i), applied to $x_2,x_3,x_4$, 
this yields that $x_4$ is uniquely
determined, and we have 
$x_4 \in S^1$, and the rotation about $o$ taking $x_1$ to $x_2$ takes $x_2$
to $x_3$ and $x_3$ to $x_4$.
Continuing like this, we conclude that
$K$ is a regular $n$-gon
inscribed into $B^2$. 
$\blacksquare $
\enddemo

\head Area sum for normalized $o$-symmetric planar convex bodies\endhead

Recall that in [BMMR] the main tool was the lower estimate for polar sectors
of convex bodies in ${\Bbb{R}}^2$. There we had a sector $xoy$, with vertex at
$o$, and of angle in
$(0, \pi )$, with $x,y \in \partial K$, and we had two supporting lines of
$K$ at $x$ and $y$, respectively, intersecting in the half-plane bounded by
${\text{aff}} \{ x,y \} $ not containing $K$.
Thus the sector contained $[x,o,y]$ and was
contained in the convex quadrangle with two sides $[o,x]$ and $[o,y]$, and
other two sides lying on the given supporting lines. This quadrangle could be
an arbitrary convex quadrangle, up to affinities. (Observe that the only affine
invariants of a convex quadrangle are the two ratios in which the two
diagonals bisect each other.) Below we will use a convex deltoid, that is up to
affinities characterized by the fact that one diagonal bisects the other one
in its midpoint. However, this will be sufficient for our theorems. In fact
we cannot solve the general case about the maximum polar area, this remains an
open 

\newpage

question.

Here we proceed analogously.
First we maximize the area product for an angular domain of angle $2 \alpha \in
(0, \pi )$, say $\angle (\cos \alpha , - \sin \alpha ) o
(\cos \alpha , \sin \alpha )$, where we suppose that the convex sector
contains the triangle $T := [(\cos \alpha , - \sin \alpha ), o, 
(\cos \alpha , \sin \alpha )]$, and is contained in the deltoid  
$Q := [(\cos \alpha , - \sin \alpha ), o, 
(\cos \alpha , \sin \alpha ), (1/\cos \alpha , 0)]$. Observe that the tangents
to $S^1$ at $(\cos \alpha , \pm \sin \alpha )$ contain the sides 
$[(\cos \alpha , \pm \sin \alpha ),$
\newline 
$(1/\cos \alpha , 0)]$ of our deltoid.
Of course, our arguments will be linearly invariant, as is the volume product,
due to the formula 
$(TK)^*=({T^{*}})^{-1}K^*$ for a {\it{non-singular}} linear transformation $T$.

For this we need Steiner symmetrization.
For an $o$-symmetric planar convex body $K$, and
a line $l$ through $o$, translate each chord of $K$
orthogonal to $l$ 
(including the ones that degenerate to points)
orthogonally to $l$ in such a way
that the midpoint of the translated chord should lie on $l$.
The union of these translates
is the Steiner symmetral $K'$ of $K$ with respect to $l$. 
Clearly $K'$ is an $o$-symmetric planar convex body with $V(K')=V(K)$.
According
to K.M. Ball's PhD thesis [...] (see also M. Meyer, A. Pajor [...]), we have
$$
V \left( (K')^* \right) \geq V(K^*),
\tag P1
$$
with equality if and only if
$K'$ is a linear image of $K$, by a linear map preserving all straight lines
orthogonal to $l$. (A similar statement holds, without $o$-symmetry, for
$V \left( \left( C-s(C) \right) ^* \right)$, with ``linear'' replaced by
``affine'', cf. ???). 




For the following proposition we need some notations.

Let $\alpha \in (0, \pi /2 )$, and let $a:= (\cos
\alpha , - \sin \alpha ), b:=(1/ \cos \alpha , 0), c:=(\cos
\alpha ,  \sin \alpha )$. Let $T:=[o,a,c]$ and $Q:=[o,a,b,c]$. Then $a,c \in
S^1$ and
${\text{aff}}\,\{ a,b \} $ and ${\text{aff}}\,\{ b,c \} $ are tangents to 
$S^1$. Let $K$ be an $o$-symmetric 
convex body with $a,c \in \partial K$, and let the counterclockwise arcs
$I:=c(-a)$ and $-I=(-c)a$ of $S^1$ be contained in $\partial K$. 
Then 
$$
[I, -I] \subset K \subset B^2 \cup [a,b,c] \cup [-a,-b,-c] \,.
$$
Then we have also 
$$
[I, -I] \subset K^* \subset B^2 \cup [a,b,c] \cup [-a,-b,-c] \,.
$$

Let $C:=K
\cap Q$ and $C^*:=K^* \cap Q$. Observe that $\partial C^* $ also contains
$I \cup (-I)$, and that $T \subset C \subset Q$.
Further, we have that $K$ or $K^*$ is the union of $C \cup (-C)$ or
$C^* \cup (-C^*)$,
and the two sectors $[ \pm I,o]$ of $B^2$.
Observe that we have $V(C) \in [\cos \alpha \sin \alpha , \tan \alpha ]$. Both
the minimal and maximal values of $V(C)$ are attained for a unique $K$: namely
for $K=[I, -I]$ and for $K=B^2 \cup [a,b,c] \cup [-a,-b,-c]$.
Still observe
$$
V(K) = 2V ( [o,I] ) + 2V(C)  {\text{ and }} 
V(K^*) =  2V ( [o,I] )  + 2V(C^*) \,.
\tag P2
$$
Hence maximization of $V(C^*)$ is equivalent to maximization of $V(K^*)$.


\proclaim{Proposition} {\it{With the above notations, let $V(C) \in 
(\cos \alpha \sin \alpha , \tan \alpha )$ be fixed. 

\newpage

Then
the maximum of $V(C^*)$ occurs, e.g., for the following cases.
\newline
(i) for $V(C) < \alpha $ e.g., for $\partial C$ 
being the union of $[o,a] \cup
[o,c]$ and an ellipsoidal arc $J$ in $B^2$ joining $a$ and $c$ (in the
positive sense), the ellipse having as centre $o$;
\newline
(ii) for $V(C) = \alpha $ e.g., for $C=B^2 \cap Q$;
\newline 
(iii) for $V(C) > \alpha $ e.g., for $\partial C$ 
being the union of $[o,a] \cup
[o,c]$ and
two segments of equal length
$[a,a'] \subset {\text{aff}}\, \{ a,b \} $ and 
$[c,c'] \subset {\text{aff}}\, \{ c,b \} $, and 
an ellipsoidal arc $J$ joining $a'$ and $c'$ (in the
positive sense), the ellipse having as centre $o$, and having ${\text{aff}}\,
\{ a,b \} $ and ${\text{aff}}\,\{ b,c \} $ as supporting lines.}}
\endproclaim


Observe that for $C$ of the form (i) or (iii) we have that
$C^*$ is of the form in (iii) or (i), respectively (for suitable areas).


\demo{Proof}
{\bf{1.}}
We begin with Steiner's symmetrization of $K$ 
with respect to the $x$ axis, obtaining $K'$. 
Observe that then $K \cap [I, -I]$ remains
invariant, and $K \cap [a,b,c]$ will be replaced by $K' \cap [a,b,c]$, which
is just the Steiner symmetral of $K \cap [a,b,c]$ with respect to the $x$
axis. We introduce the notation 
$$
C' := [o, \left( K' \cap [a,b,c] \right) ] \,.
$$
Then 
$$
\cases
V(K)=2V \left( [o,I] \right)  + 2V(C) =
V(K')=2V \left( [o,I] \right)  + 2V(C')
{\text{ and}}\\
V(K^*)=2V \left( [o,I] \right)  + 2V(C^*) =
V \left( (K')^* \right) = 2V \left( [o,I] \right)  + 2V \left( (C') ^*
\right) \,.
\endcases
\tag P3
$$
Hence maximization of $V(C^*)$ is equivalent to maximization of $V(K^*)$.
Still observe that by the equality case of 
\thetag{P1} $V(K^*)$ can attain its maximum only in the
case when $K'=K$, since the only linear map preserving vertical lines and
taking $[I, -I]$ to itself is the identity.

Therefore we may suppose that $K$ and thus also $C$ is symmetric with respect
to the $x$ axis. In particular, 
$$
\cases
[a,b] \cap \partial K  {\text{ and }} [b,c]
\cap \partial K  {\text{ will be symmetric images of each}} \\ 
{\text{other with respect to the }} x 
{\text{ axis, thus they have the same length}}\,.
\endcases
\tag P4
$$
This length can be $0$ or positive, but in the second case
by the hypothesis
$V(C) \in (\cos \alpha \sin \alpha , \tan \alpha )$ of the proposition 
its length
is smaller than the length of $[a,b]$. 

{\bf{2.}}
Now let us consider a convex body $K$ satisfying the hypotheses of the
Proposition for which $V(K^*)$ is maximal, and the body
$C$ corresponding to this $K$. Then by {\bf{1}}
both $K$ and $C$ are necessarily symmetric with
respect to the $x$-axis.
We have two cases.
\newline
(i) $[a,b] \cap \partial K$
(and then also $[b,c] \cap \partial K$) has a positive length.
\newline
(ii)
$[a,b] \cap \partial K$
(and then also $[b,c] \cap \partial K$) is a point.

For this $K$ and $C$ let us inscribe to 
$[(\partial C) \setminus T] \cup \{ a, c \} $ a convex 
polygonal arc
$x_1 \dots x_n$, 
in positive orientation, in the following way. We have
$x_1=a$ and $x_n=c$. In case (i) the segment
$[a,b] \cap \partial K$ (and $[b,c] \cap \partial K$)
should be $[x_1,x_2]$ (or $[x_{n-1}x_n]$, respectively). The further points
$x_3, \dots , x_{n-2}$ in case (i) or $x_2, \dots , x_{n-1}$ in case (ii)
are placed so that they divide the counterclockwise
arc $\widehat{x_2 \dots x_{n-1}}$ or 
$\widehat{x_1 \dots x_n}$ 

\newpage

of $\partial C$ to subarcs of equal length.
This guarantees symmetry of our
polygonal arc with respect to the $x$ axis. 
Since the length of the arc $[\partial C \setminus T] \cup \{ a,c \} $ 
is at most
the length of the arc ${\widehat{abc}}$, 
therefore the length of the subarcs ${\widehat{x_ix_{i+1}}}$,
excepting ${\widehat{x_1x_2}}$ and ${\widehat{x_{n-1}x_n}}$ in case (i), is
$O(1/n)$.
 
We write
$$
\cases
K_n := [I, -I] \cup {\text{conv}} \{ x_1, \dots , x_n \} \cup
{\text{conv}} \{ -x_1, \dots , -x_n \} \\
{\text{and }}
C_n:= [o, K_n \cap [a,b,c] ] \,.
\endcases
$$
Then $K_n \to K$ and $C_n \to C$ in the Hausdorff metric in both cases (i)
and (ii). In fact, we have on one hand $K_n \subset K$ and $C_n \subset C$.
On the other hand, if some $p \in C$ is separated from $C_n$ by a side
$x_ix_{i+1}$ (in case (i) this cannot be $x_1x_2$ or $x_{n-1}x_n$), 
then 
$$
\min \{ |x_ip|, |x_{i+1}p| \} \le (|x_ip| + |x_{i+1}p|)/2 \le
|{\widehat{x_ix_{i+1}}}| / 2 = O(1/n)\,.
$$ 
Hence the Hausdorff distance of $K$ and $K_n$, as well as that of $C$ and
$C_n$ is $O(1/n)$.

Now we define ${\tilde K}_n := [I, -I, J, -J]$, where $J$ is 
a strictly convex
polygonal arc $J:= {\tilde x}_1, \dots , {\tilde x}_k$ 
which has the following properties. We have $k \le n$.
Further, $J$
lies in $[a,b,c]$ and satisfies $x_1=a$ and
$x_k=c$, and in case (i) still $[x_1,x_2] \supset [a,b] \cap \partial K$ and 
$[x_{k-1},x_k] \supset [b,c] \cap \partial K$, 
and $V \left( [ {\tilde x}_1, \dots , {\tilde x}_k ] \right) =
V \left( [ x_1, \dots , x_n ] \right) $. Lastly, $J$ is such a one among all 
polygonal arcs satisfying all above listed properties, for which
$V\left( ({\tilde K}_n)^* \right)$ attains its maximal value.
(By a {\it{(strictly) convex polygonal arc}}
we mean one for which the sides following each other turn (strictly) 
counterclockwise.)
Observe that till now we could not exclude $[{\tilde{x_1}},{\tilde{x_2}}] \ne 
[a,b] \cap \partial K$ and 
$[{\tilde{x}}_{k-1},{\tilde{x}}_k] \ne 
[b,c] \cap \partial K$ 
in either case (i) or case (ii), and also we do not yet know $k=n$. 

Then $V\left( (K_n)^* \right) 
\le V\left( ({\tilde K}_n)^* \right) $. Then for some subsequence 
${\tilde{K}}_{n(i)}$ of the ${\tilde{K}}_n$'s there exists
a limit convex body, and also the equal lengths of the
intersections ${\tilde{K}}_{n(i)} \cap [a,b]$ and ${\tilde{K}}_{n(i)} 
\cap [b,c]$ are convergent.
Then also $V \left( (\lim {\tilde K}_{n(i)})^{^*} \right) $ 
will have a maximal value among the considered convex bodies $K$.

{\bf{3.}}
Now we want to apply the method of proof of Theorem A. Its part {\bf{1}}
consisted of three
main steps. First we showed $k=n$. Then we showed (i) and (ii) from the proof
of Theorem A.

Now, rather than (i) and (ii) from {\bf{2}} we have to distinguish between the
following cases.
\newline
(i') ${\tilde{x}}_2 \in (a,b)$ and ${\tilde{x}}_{k-1} \in (b,c)$, and
\newline
(ii') ${\tilde{x}}_2 \not\in (a,b)$ and ${\tilde{x}}_{k-1} \not\in (b,c)$.
Observe that by strict convexity, in case (i') we have
$$
{\tilde{x}}_3, \dots , 
{\tilde{x}}_{n-2} \in {\text{int}}\,[a,b,c]\,,
$$
while in case (ii') we have
$$
{\tilde{x}}_2, \dots , 
{\tilde{x}}_{n-1} \in {\text{int}}\,[a,b,c]\,.
$$

As soon as some vertex ${\tilde{x_i}}$ lies in 
${\text{int}}\,[a,b,c]$, it is freely
movable till some small distance. Therefore the statements in the proof of
Theorem A, {\bf{1}}, which used small movements of such $x_i$'s, remain valid
also here.

\newpage

In particular, (i) of Theorem A, {\bf{1}} remains true for any three
consecutive vertices ${\tilde{x}}_i,{\tilde{x}}_{i+1},{\tilde{x}}_{i+2}$, for
$2 \le i \le k-3$ in case (i'), and for $1 \le i \le k-2$ for case (ii').

However, for the proof of $k=n$ and (ii) in the proof of
Theorem A, {\bf{1}} the area of $K$ was not preserved, so we have to take more
care in the proof here.

{\bf{4.}}
We begin with the proof of $k=n$. In the proof of Theorem A, {\bf{1}} we cut
off from $K$
a part of area $\Theta (\varepsilon ^2)$, so for keeping the area constant,
we have to give this area back at some other place.

Namely, supposing $k<n$ consider a straight line $l$ such that ${\tilde K}_n
\cap l = \{ {\tilde{x}}_{i+1} \} $. 
Then push parallelly $l$ towards $o$ through a small
distance $\varepsilon =\varepsilon ({\tilde K}_n) > 0$, obtaining $l'$. 
Then consider the
$o$-symmetric convex body 
$({\tilde K}_n)'\,\,( \subset {\tilde K}_n)$ obtained from ${\tilde K}_n$ by
cutting ${\tilde x}_{i+1}$ and $-{\tilde x}_{i+1}$
from ${\tilde K}_n$ by $l'$ and by $-l'$. Then $({\tilde K}_n)'$ has in the
closed angular domain $\angle aoc$ $k+1 \le n$ vertices, and 
$V({\tilde K}_n)-
V\left( ({\tilde K}_n)'\right) = \Theta (\varepsilon ^2)$. 
Moreover, $({\tilde K}_n)'$ has vertices
${\tilde x}_i,{\tilde x}_{i+1}',{\tilde x}_{i+1}'',{\tilde x}_{i+2}$ 
in positive sense. 
Now fix ${\tilde x}_i$ and aff\,$\{ {\tilde x}_{i+1}',{\tilde x}_{i+2}''
\} $, and rotate the side line aff\,$\{ {\tilde x}_i,{\tilde x}_{i+1}' \} $ 
about ${\tilde x}_i$
outwards through an
angle $\Theta (\varepsilon ^2)$ (and making the analogous change also at
$-{\tilde x}_i$) so that for the new $o$-symmetric convex body
$({\tilde K}_n)''$, that has also $k+1$ vertices in the closed
angular domain $\angle aoc$
we should have $({\tilde K}_n)'' \supset ({\tilde K}_n)'$ and 
$V \left( ({\tilde K}_n)'' \right) =V({\tilde K}_n )$. 
Now we turn to the polar bodies. By ${\tilde K}_n \supset ({\tilde K}_n)' 
\subset ({\tilde K}_n)''$ we have 
$({\tilde K}_n)^* \subset \left( ({\tilde K}_n)' \right) ^* \supset 
\left( ({\tilde K}_n)'' \right) ^*$. Then $({\tilde K}_n)^*$ 
has neighbouring vertices
$({\tilde x}_i{\tilde x}_{i+1})^*$ and $({\tilde x}_{i+1}{\tilde x}_{i+2})^*$, 
connected by the side on the line ${\tilde x}_{i+1}^*$. Then
$\left( ({\tilde K}_n)' \right) ^*$ will have a new vertex $(l')^*$, 
at a distance $\Theta (\varepsilon )$ 
from ${\tilde x}_{i+1}^*$, hence 
$V\left( \left( ({\tilde K}_n)' \right) ^* \right) - V\left( ({\tilde K}_n)^*
\right) = \Theta (\varepsilon )$.
The rotation of the side line aff\,$\{ {\tilde x}_i,{\tilde x}_{i+1}' \} $ 
about ${\tilde x}_i$ through
an angle $O (\varepsilon ^2)$ implies motion of $($aff\,$\{ {\tilde x}_i,
{\tilde x}_{i+1} \} )^*$ 
on the line $({\tilde x}_i)^*$ 
through a distance $O (\varepsilon ^2)$, hence 
$V \left( \left( ({\tilde K}_n)' \right) ^* \right) 
-V \left( \left( ({\tilde K}_n)'' \right) ^* \right) = O (\varepsilon ^2)$.
Therefore,
$$
V \left( \left( ({\tilde K}_n)'' \right) ^* \right) 
= V \left( ({\tilde K}_n)^* \right) + 
\Theta (\varepsilon ) - O (\varepsilon ^2) > V \left( ({\tilde K}_n)^*
\right),
$$
for $\varepsilon >0$ sufficiently small, a contradiction. This contradiction
proves our claim $k=n$.

{\bf{5.}}
For the proof of (ii) in the proof of
Theorem A, {\bf{1}} (by dualizing (i) in the proof of
Theorem A), the area of $K$ was not preserved, even just conversely, the area of
the polar was preserved. Therefore we have to give a new proof, which
preserves the area of $K$.

So we consider four consecutive vertices ${\tilde x}_i,{\tilde x}_{i+1},
{\tilde x}_{i+2},{\tilde x}_{i+3}$ of ${\tilde K}_n$, and $m$ is the
intersection point of ${\text{aff}}\,\{ {\tilde x}_i,{\tilde x}_{i+1} \} $
and ${\text{aff}}\,\{ {\tilde x}_{i+2},{\tilde x}_{i+3} \} $. Then $m$ is a
finite point, separated from ${\tilde K}_n$ by 
${\text{aff}}\,\{ {\tilde x}_{i+1},{\tilde x}_{i+2} \} $.

We may suppose that the lines ${\text{aff}}\,\{ {\tilde x}_i,{\tilde x}_{i+1}
\} $, ${\text{aff}}\,\{ {\tilde x}_{i+1},{\tilde x}_{i+2} \} $ and 
${\text{aff}}\,\{ {\tilde x}_{i+2},$
\newline
${\tilde x}_{i+3} \} $ are tangent to $S^1$,
and also that $m$ lies on the positive $x$ axis. So ${\text{aff}}\,\{ o, m \}
$ is the $x$ axis.
Then ${\text{aff}}\,\{ {\tilde x}_i,{\tilde x}_{i+1} \} $ and 
${\text{aff}}\,\{ {\tilde x}_{i+2},{\tilde x}_{i+3} \} $ are symmetric with
respect to 

\newpage

the $x$ axis. Then of course, also their polars
$({\text{aff}}\,\{ {\tilde x}_i,{\tilde x}_{i+1} \})^*$ and 
$({\text{aff}}\,\{ {\tilde x}_{i+2},{\tilde x}_{i+3} \} ) $ are symmetric with
respect to the $x$ axis, thus 
$$
\cases
{\text{the line connecting }} 
({\text{aff}}\,\{ {\tilde x}_i,{\tilde x}_{i+1} \})^*  \\
{\text{and }} ({\text{aff}}\,\{ {\tilde x}_{i+2},{\tilde x}_{i+3} \} )^* 
{\text{ is vertical}}\,.
\endcases
\tag P5
$$
We may assume, for contradiction, that $({\tilde x}_
{i+1} + {\tilde x}_{i+2})/2$ lies below the $x$ axis.

We begin with rotating the side line 
${\text{aff}}\,\{ {\tilde x}_{i+1},{\tilde x}_{i+2} \} $ about the midpoint 
$({\tilde x}_{i+1} + {\tilde x}_{i+2})/2$ of the side $[{\tilde x}_{i+1},
{\tilde x}_{i+2}]$, through some small angle $\varepsilon >0$ 
in positive sense. 
(That is, we move toward the symmetric position, like in the proof of (i) in
{\bf{1}} of the proof of Theorem A.) The new positions of ${\tilde x}_{i+1}$
and ${\tilde x}_{i+2}$ are denoted by $({\tilde x}_{i+1})'$ and
$({\tilde x}_{i+2})'$.
The body obtained by this change is
denoted by $({\tilde K}_n)'$. Then we have $V \left( ({\tilde K}_n)' \right)
- V ({\tilde K}_n) = O( \varepsilon ^2)$. 

Of course we have to still to restore the original area. This is done by a
parallel replacement of the already rotated side line through a distance
$O( \varepsilon ^2)$. 
The new positions of $({\tilde x}_{i+1})'$ and 
$({\tilde x}_{i+2})'$ will be denoted 
by $({\tilde x}_{i+1})''$ and $({\tilde x}_{i+2})''$, and the
thus obtained body will be denoted by
$({\tilde K}_n)''$.

Observe that the centre of rotation $({\tilde x}_{i+1} + {\tilde x}_{i+2})/2$
of our first motion lies in $[m,s,t] \setminus B^2 \subset [m,s,t]$, 
where $s$ and $t$ are the 
points of tangency of $S^1$ with the side lines 
${\text{aff}}\,\{ {\tilde x}_i,{\tilde x}_{i+1} \} $ and 
${\text{aff}}\,\{ {\tilde x}_{i+2},{\tilde x}_{i+3} \} $.
Therefore $({\tilde x}_{i+1} + {\tilde x}_{i+2})/2$ lies in the open right
halfplane given by $x>0$. By hypothesis, it also lies in the open lower
half-plane given by $y<0$. So it lies in the open fourth coordinate quadrant.
Therefore 
$$
{\text{the slope of the polar line }}
[({\tilde x}_{i+1} + {\tilde x}_{i+2})/2]^* {\text{ lies in }} (0, \infty )\,.
\tag P6
$$

Now we consider the polars. We begin with the first part of the motion, i.e.,
with the rotation of the side line about the side midpoint.
Rotation of the side line
${\text{aff}}\,\{ {\tilde x}_{i+1},{\tilde x}_{i+2} \} $ about the midpoint 
$({\tilde x}_{i+1} + {\tilde x}_{i+2})/2$ of the side $[{\tilde x}_{i+1},
{\tilde x}_{i+2}]$ implies moving the polar point 
$({\text{aff}}\,\{ {\tilde x}_{i+1},
{\tilde x}_{i+2} \} )^*$ on the polar line
$[({\tilde x}_{i+1} + {\tilde x}_{i+2})/2]^*$, so that this point moves
counterclockwise, when looked upon from $o$. Therefore by \thetag{P5}
its $x$ coordinate 
increases, by a quantity $\Theta ( \varepsilon )$. 

The polar body $({\tilde K}_n)'$ has consecutive vertices  
$({\text{aff}}\,\{ {\tilde x}_i,{\tilde x}_{i+1} \})^*$,
$\left( {\text{aff}}\,\{ ({\tilde x}_{i+1})',\right.$
\newline
$\left. ({\tilde x}_{i+2})' \} \right) ^*$
and $({\text{aff}}\,\{ {\tilde x}_{i+2},{\tilde x}_{i+3} \} ) $. The first and
third of these vertices are fixed, and their connecting line is vertical, by
\thetag{P5}. The second vertex lies on the right hand side of this vertical
line, and in its new position 
$\left( ({\text{aff}}\,\{ {\tilde x}_{i+1},{\tilde x}_{i+2} \} )' \right) ^*$
its $x$ coordinate is greater than in its original position 
$ ({\text{aff}}\,\{ {\tilde x}_{i+1},{\tilde x}_{i+2} \} ) ^*$, namely the
difference is $\Theta ( \varepsilon )$. Therefore 
$$
V\left (K') ^* \right) - V (K ^*) =  \Theta ( \varepsilon ) \,.
\tag P7
$$

Now we consider the second part of this motion, i.e., the parallel
displacement of the already rotated side, through a distance $O( \varepsilon
^2)$. Then the vertices
$({\text{aff}}\,\{ {\tilde x}_i,{\tilde x}_{i+1} \})^*$,
and $({\text{aff}}\,\{ {\tilde x}_{i+2},{\tilde x}_{i+3} \} )^* $ remain fixed,
but the vertex 
$\left( {\text{aff}}\,\{ ({\tilde x}_{i+1})',\right.$

\newpage

$\left. ({\tilde x}_{i+2})' \} \right) ^*$
moves to its new position 
$\left( {\text{aff}}\,\{ ({\tilde x}_{i+1})'',({\tilde x}_{i+2})'' \} 
\right) ^*$, and the distance of these last two points is $O( \varepsilon
^2)$. Hence 
$$
V\left( (K'') ^* \right) - V\left( (K') ^* \right) =O( \varepsilon ^2)
\tag P8
$$
Then \thetag{P7} and \thetag{P8} imply
$$
V\left( (K'') ^* \right) = V ( K ^* ) + \Theta (\varepsilon ) + O( \varepsilon
^2 ) = V ( K ^* ) + \Theta (\varepsilon ) > V ( K ^* )\,,
\tag P9
$$
that is a contradiction. This proves (ii) 
in the proof of Theorem A, {\bf{1}}.

{\bf{8.}}
By {\bf{5}}, {\bf{6}} and {\bf{7}} we have that the extremal polygonal line
$x_1 \dots x_n$ has exactly $n$ vertices, 
and both (i) and (ii) from the proof of Theorem
A, {\bf{1}} hold.

Recall the cases (i') and (ii') from {\bf{3}}. 
From {\bf{3}}, in case (i') we have
$$
{\tilde{x}}_3, \dots , 
{\tilde{x}}_{n-2} \in {\text{int}}\,[a,b,c]\,,
$$
while in case (ii') we have
$$
{\tilde{x}}_2, \dots , 
{\tilde{x}}_{n-1} \in {\text{int}}\,[a,b,c]\,.
$$

Then, like in the proof of Theorem A, we obtain that in case (i') 
${\tilde x}_2 \dots
{\tilde x}_{n-1}$, while in case (ii) ${\tilde x}_1 \dots
{\tilde x}_n$ is
inscribed to an elliptical arc, where the respective
ellipse $E_n$ has centre $o$, and, in the
metric of $E_n$, the rotation carrying $x_2$ to $x_3$ in case (i), or carrying
$x_1$ to $x_2$ in case (ii) also carries all
$x_i$ to $x_{i+1}$, for $3 \le i \le n-2$, or $2 \le i \le n-2$, respectively.
Since any two metrics on
${\Bbb{R}}^2$ are equivalent, therefore changing the original metric of
${\Bbb{R}}^2$  
to the metric determined by $E_n$ preserves lengths and areas,
up to an at most constant positive factor. Therefore we may use henceforward
for fixed $n$ the metric of $E_n$.

Now let us replace the polygonal arc ${\tilde x}_3 \dots {\tilde x}_{n-2}$ in
case (i'), and ${\tilde x}_2 \dots {\tilde x}_{n-1}$ in case (ii') by the
respective arc of $E_n$ (and similarly for their mirror images with respect to
$o$). Then in both cases we obtain a convex body ${\tilde L}_n$, 
since the chords 
$[{\tilde x}_2 , {\tilde x}_3]$ and $[{\tilde x}_{n-2} , {\tilde x}_{n-1}]$
in case (i'), and 
$[{\tilde x}_1 , {\tilde x}_2]$ and $[{\tilde x}_{n-1} , {\tilde x}_n]$
in case (ii') span lines intersecting in
$[a,b,c]$.

Then we have ${\tilde K}_n \subset
L_n$, and the Hausdorff distance of ${\tilde K}_n$ and $L_n$ 
is $O (1/n^2) = o(1)$. For a
suitable subsequence $n(i)$ (cf. {\bf{2}}) we have that 
\newline
(a) ${\tilde K}_n$ is
convergent, and then necessarily ${\tilde L}_{n(i)}$ has the same limit, and
\newline
we have case (i') for all $n(i)$ or case (ii') for all $n(i)$,
and 
\newline
(b) in case (i') the points ${\tilde x}_2 \in (a,b)$ and ${\tilde x}_{n-1} \in
(b,c)$ are convergent, and also 
\newline
(c) $E_{n(i)}$ is convergent. 
\newline
Therefore we claim
that we may use
the metric of any $E_n$, or the metric of $E:= \lim _{i \to \infty }
E_{n(i)}$, the changes in lengths and areas are still bounded by an at most
constant factor. 

For this we have to show that all $E_n$'s lie in a compact family. All of them 
pass through all four vertices $( \pm \cos \alpha , \pm \sin \alpha $
of an axis-parallel rectangle, hence have
axisparallel axes themselves. The positive horizontal
semiaxes have endpoints in ${\text{int}}\,[a,b,c]$,
hence are bounded from below and from above. The positive vertical

\newpage

semiaxes are bounded from below by the $y$-coordinate of the vertex ${\tilde
x}_1$ in case (i') and of the vertex ${\tilde x}_2$ in case (ii'). In case
(i') this is fixed. Also in case (ii') these
cannot be arbitrarily small, since then the areas of $[{\tilde x}_1, \dots
,{\tilde x}_n]$ would tend to $V([a,b,c])$, which was excluded. 
Also they cannot be arbitrarily large. Namely then the chords of $E_n$ from
the endpoints of their vertical semiaxes with positive (negative)
$y$-coordinates to the endpoints with positive $x$-coordinates 
and positive (negative) $y$-coordinates of the
above considered elliptical arcs span lines on whose left hand side lies
$[{\tilde x}_1, \dots , {\tilde x}_n]$. Then the areas of 
$[{\tilde x}_1, \dots , {\tilde x}_n]$ would tend to $0$, which 
was also excluded.

Therefore the limit of the ${\tilde L}_{n(i)}$'s exists. It is $o$-symmetric.
In the closed
angular domains $[ \alpha, \pi - \alpha ]$ and its mirror image with
respect to $o$ it is bounded by the respective arcs of $S^1$. In the closed
angular
domains $[- \alpha , \alpha ]$ it is bounded in case (i') by the segment
$[\lim {\tilde x}_1, \lim {\tilde x}_2]$, its mirror image w.r.t. the
$x$-axis, and an elliptical 
arc symmetric w.r.t. the $x$-axis, while in case (ii')
only with the hyperbola arc symmetric w.r.t. the $x$-axis. In case (i')
possibly $\lim {\tilde x}_1 = \lim {\tilde x}_2$ --- this case we count to
case (ii').

In case (ii') the area of the convex hull of the arc of ellipse (that is
equal to $V([{\tilde x}_1, \dots , {\tilde x}_n ])$ for each $n$) uniquely
determines the arc of ellipse. Namely the boundary 
of the ellipse passes through the four points $(\pm \cos \alpha , \pm
\sin \alpha )$. Therefore any two of them does not have any further
intersection points, so their parts in the open angular domain $( -\alpha ,
\alpha )$ are disjoint, hence they are linearly ordered by inclusion. The arc
giving the minimal area is
$[a,c]$, which however cannot be attained by the restriction on the
area. We assert that 
the maximal of them is the arc of $S^1$ in this angular domain. Namely
any two such ellipses have a transversal intersection point at $a$ (and $c$). 
Else they would have eight common points with multiplicities, hence would be
equal. In particular, the arc of $S^1$ in this angular domain and any
elliptical arc passing through the four points $(\pm \cos \alpha , \pm
\sin \alpha )$ have different tangents at $a$. In case of strict inclusion of
the sector of $B^2$ and the respective sector of the elliptical arc the
tangent of the elliptical arc at $a$ would point outside from the quadrangle
$Q$, that is impossible, since it lies in $Q$.

In case (i') the situation is more complicated. Then we have on the boundary
of ${\tilde L}$ segments of equal positive lengths 
$[{\tilde x}_1, {\tilde x}_2]$ and $[{\tilde x}_{n-1}, {\tilde x}_n]$. On
$\partial L$, at 
${\tilde x}_2$ and its mirror image w.r.t. the $x$-axis 
there joins to these segments
an arc of an ellipse, also symmetric w.r.t. the $x$-axis. Now we
have two parameters: the positive length of the segment $[{\tilde x}_1,
{\tilde x}_2]$, and one more parameter that distinguishes the elliptical arc
among all those elliptical arcs that pass through the four points $ \pm
{\tilde x}_2$ and their symmetric images w.r.t. the $x$-axis. 

However, here is still one restriction. The body $L$ cannot have a non-smooth
point at ${\tilde x}_2$. This can be proved similarly as formerly
already several times. Namely, we
put a line $l$ that is a supporting line of $L$ at  ${\tilde x}_2$, but is
not equal to any of the half-tangents at ${\tilde x}_2$. Now let $\varepsilon
>0$ be sufficiently small, and push $l$ inward to $L$ (and similarly at all
three symmetric images of ${\tilde x}_2$). The loss of area is $\Theta (
\varepsilon ^2)$. Since the area of $K$ must remain constant, we have to give
this loss of area somewhere back. For this consider the elliptical arc in
question, and close to the midpoint of this elliptical arc we choose a point
$p$ on the positive $x$-axis outside of $L$, having a distance $\delta $ from
$L$. We consider the convex hull $[L, p -p]$. This has area 

\newpage

$V(L) + \Theta
(\delta ^{3/2})$. Similarly, $[L, p -p]^*$ has an area $V(L^*) + \Theta
(\delta ^{3/2})$. Then choosing $\delta $ suitably, the added area will be
exactly equal to the lost area, so the area of $K$ will not change by
performing these two operations. At the same time, $\delta ^{3/2} = \Theta (
\varepsilon ^2)$. By cutting off $K$ with $l$ (and its three symmetric images)
$K^*$ will change to the convex hull of the original $K^*$ and $l^*$ and its
three symetric images. However, this means an increase in the area by $\Theta
(\varepsilon )$. Then cutting the new $K^*$ by the line $p^*$ 
(and its three symmetric images) causes the area to decrease with $\Theta 
(\delta ^{3/2}) = \Theta (\varepsilon ^2)$. That is, the original $V(K^*)$
changed to $V(K^*) + \Theta ( \varepsilon ) + \Theta ( \varepsilon ^2)=
V(K^*) + \Theta ( \varepsilon ) > V(K^*)$ a contradiction.

This contradiction implies that the extremal body $L$ is smooth at 
${\tilde x}_2$ (and its three symmetric images). This means that it is already
uiniquely determined (by passing through the point ${\tilde x}_2$ and having
given tangents there --- namely two such ellipses would have eight common
points with multiplicities, which is a contradiction). Observe that the polars
of such arcs, consisting of two symmetrically placed segments and an arc of a
hyperbola, are just elliptical arcs of the form in case (ii'),
hence the corresponding (polar) sectors are totally ordered by
inclusion. Therefore, for a given area of the sector of $K$ in the closed
angular
domain $aob$ the part of the boundary in this angular domain is uniquely
determined. By the upper 
inequality on the area of the sector of $K$ in this angular
domain the degenerate position, when the sector of $K$ would coincide with the
deltoid $[o,a,b,c]$, cannot be attained. 

We obtain case (ii) of the Proposition when both (i') and (ii') hold.
$\blacksquare $
\enddemo

The {\it{John ellipse}} $E_i(K)$, or the {\it{L\"owner ellipse}} $E_o(K)$ 
{\it{of an $o$-symmetric convex body}} $K \subset {\Bbb{R}}^d$
is a the $o$-symmetric ellipsoid of maximal volume contained in $K$, or of
minimal volume containing $K$, respectively.
Both of them are unique, and the polar of
the John ellipsoid of $K$ is the L\"owner ellipsoid of $K^*$ (using duality).
By [Be], $B^2$ is the
John ellipse or the L\"owner ellipse of $K \subset {\Bbb{R}}^2$ 
if and only if it is contained in
$K$, or contains $K$, and $(\partial K) \cap S^1$ contains the vertices of a
square with vertices $p_1,...,p_4$ in positive cyclic order, inscribed to
$S^1$, or the vertices of an $o$-symmetric convex hexagon with vertices
$p_1,...,p_6$, in positive cyclic order, where $\angle p_iop_{i+1} < \pi /2 $
(Behrend, [Behrend]).


\proclaim{Theorem B} If $B^2$ is the John ellipse
or the L\"owner ellipse of
a planar o-sym\-met\-ric convex body $K$, then 
$$
V(K)+V(K^*)\leq 2\pi.
$$
\endproclaim


\definition{Remark} This inequality is specific to $o$-symmetric
planar convex bodies. If either $K$ is a regular triangle inscribed into
$S^1$,
or 
or a regular cross-polytope circumscribed about $S^{d-1}$ (for the
John ellipsoid) or
a cube inscribed to $S^{d-1}$ (for the L\"owner ellipsoid) 
of dimension $d \ge 3$, the analogous
statement (i.e., when $\pi $ is replaced by the volume of the unit ball $B^d$
in ${\Bbb R}^d$) does not hold. 
\enddefinition

In fact, for $K$ a regular triangle inscribed to $S^1$ we have
$V(K)+V(K^*)=15 \sqrt{3}/4 = 6.4951... > 2 \pi $. Further,
observe that the John
ellipsoid of a regular cross-polytope $D^d$, circumscribed about
$B^d$, is $B^d$ (by the criterion of John [John]). 
Then, 

\newpage

by polarity, the
L\"owner ellipsoid of a cube $C^d$, inscribed to $B^d$, is $B^d$.
For $d=3$, we have $V(C^3)+V(D^3) =44 \sqrt{3} /9 = 
8.4678... > 2 V(B^3)= 8.3775... $. For $d=4,\,\,5$ already $V(D^4)=10.6666...
> 9.8696...=2V(B^4)$ and $V(D^5)=14.9071... > 10.5275...=2V(B^5)$, so
$V(D^d) > 2V(B^d)$. Further, using the recursion formulas, one easily sees that
$V(D^d)/\left( 2V(B^d) \right) $ is increasing, 
separately for even, and for odd $d \ge 4$, that proves our claim. (For
comparing the cases of dimensions $d$ and $d+2$, the increasing property is
equivalent to $(1+2/d)^{d/2} \cdot [(d+2)/(d+1)] \cdot (2/ \pi ) > 1$, 
which holds even omitting the second factor.)


\demo{Proof of Theorem B} 
We may assume that $B^2$ is the L\"owner ellipse of $K$.
It is sufficient to prove that
if $p:=p_i,q:=p_{i+1} \in S^1 \cap \partial K$,
and the angle of the vectors $p$ and $q$ is
$2 \alpha \in (0, \pi /2]$, then
$$
V(K\cap S)+V(K^*\cap S)\leq 2 \alpha
\tag *
$$
for the convex cone $S:=\{ tp+sq :\,t,s \geq 0 \} $. 
Namely, summing all four, or all six such inequalities, we obtain the
statement of the theorem.

We may suppose $p=(1,0)$.
Let $r := (0,1) \in S^1$, 
and let $\widetilde{S}=\{ tr+sq:\,t,s \geq 0 \} $. For 
$$
\widetilde{K}_{pq}=[\pm(S\cap K) , \pm(\widetilde{S} \cap B^2)],
$$
we have
$$
{\widetilde{K}}_{pq}^*=[\pm (S \cap K^*) ,
\pm({\widetilde{S}} \cap B^2) ,  \pm \sqrt{2}\,(r-p)  ].
$$
Therefore (*) is equivalent to
$$
V({\widetilde{K}}_{pq})+V({\widetilde{K}}_{pq}^*) \leq 3+\pi \,.
$$
Let $E$ be an $o$-symmetric ellipse such that
$p,r \in \partial E$, and for the part $M$ between
the chords  $[-r,p]$ and $[-p,r]$, we have
$V(M)=V({\widetilde{K}}_{pq}) \le 1+ \pi /2$. By $\pm p,\pm r \in \partial E$
we have $V(E) \ge \pi $.
It follows from the Proposition that
$V({\widetilde{K}}_{pq}^*) \leq V(M^*)$,
and hence (*) follows from
$$
V(M \cap S_0)+V(M^* \cap S_0) \leq \pi/2 \,,
\tag **
$$
where $S_0:=S \cup {\widetilde{S}}=\{ tp+sr:\,t,s \geq 0 \} $ and $V(A \cap S_0)
\le \pi /4$. Thus it suffices to show $(**)$, which we are going to do. 

Let $a \le 1$ and $b \ge 1$ be the half axes of $E$. For $E=B^2$ $(**)$ is
fulfilled, therefore we may assume $a < 1 < b$. 
We choose a new orthonormal system of coordinates 
with basic vectors $(1/ {\sqrt{2}}, 1/ {\sqrt{2}}), 
(-1/ {\sqrt{2}}, 1/ {\sqrt{2}})$, in this order.
Then $E=\Phi B^2$, for a diagonal matrix 
$$
\Phi = \left( 
\matrix
a & 0 \\
0 & b 
\endmatrix \right) .
$$
In particular,
$E \cap S_0=\Phi \sigma $,
where $\sigma \subset S_0$ is a sector of $B^2$, having an acute angle 
$2 \alpha $ 

\newpage

at $0$, and being 
symmetric with respect to the perpendicular bisector of $[p,r]$.
(For $ 2\alpha < \pi /2$ observe $V(E) > \pi $, and hence $V(E \cap S_0) \le
\pi /4 < V(E)/4$.)
It follows that 
$\Phi \left( \cos \alpha , \sin  \alpha \right) = (a \cos \alpha , b \sin
\alpha )= (1/{\sqrt{2}},
1/{\sqrt{2}})$, hence $a=1/({\sqrt{2}}\cos \alpha )$ and $b=
1/({\sqrt{2}} \sin \alpha )$, hence
$\tan (\alpha )=a/b$ and
$$
M^*\cap S_0=[p,r,\tau ],
$$
where $\tau $ is the sector of the polar ellipse $E^*$ corresponding to the
sector $\sigma $ of $E$ by polarity.
Hence
$$
V(M \cap S_0)+V(M^* \cap S_0)=ab\alpha +
\frac{\alpha }{ab}+\frac{1- \tan ^2 \alpha }{1+\tan ^2  \alpha } .
$$
In particular, $V(M \cap S_0)+V(M^*\cap S_0)=f( \alpha )$ 
for $\alpha \in (0, \pi /4]$, where $f ( \pi /4)= \pi /2$, 
and
$$
f( \alpha )=\alpha \left( \sin (2 \alpha ) + 1 / \sin (2 \alpha ) \right) 
+ \cos (2 \alpha ) .
$$

We write $\beta := 2 \alpha $ and $g( \beta ):=f ( \beta /2 )$, where $\beta \in
(0, \pi /2]$.
The inequality $g'( \beta )>0$ becomes, rearranging, $\beta < \tan
\beta $, for $\beta \in (0, \pi /2)$, hence it holds.
Therefore
$$
V(M\cap S_0)+V(M^* \cap S_0)=g(\beta ) \le  g( \pi /2 )=\pi/2,
$$
completing the proof of (**), and in turn
the proof of Theorem~B. 
$\blacksquare $
\enddemo

\head Stability of the Blaschke-Santal\'o inequality\endhead

\proclaim{Theorem D} Let $K$ be a planar convex body, and $s(K)$ the
Santal\'o point of $K$.
If the Banach-Mazur
distance of $K$ from the ellipses is
at least $1+\varepsilon $, where $\varepsilon > 0$, then 
$$
V(K)V \left( K - s(K) \right) ^*) \leq \pi ^2 (1-c \varepsilon ^6),
$$
where $c$ is a positive absolute constant. If additionally $K$ is
$0$-symmetric, we have 
$$
V(K)V \left( K - s(K) \right) ^*) = V(K)V(K^*) \leq 
\pi ^2 \left( 1-\left( 8 + o(1) \right) \varepsilon ^3 \right),
$$
for $\varepsilon \to 0$.
\endproclaim

\demo{Proof} 
It follows by \cite{B\"o}, Theorem 1.4,
together with \cite{MR06}, Theorem 1 (and Theorem 13),
that there exists an $o$-symmetric planar
convex body $\widetilde{K}$
with $V(\widetilde{K})V(\widetilde{K}^*)$
\newline
$\geq V(K)V \left( \left( K - s(K) \right) ^* \right) $,
whose Banach-Mazur distance from the ellipses is
at least $1+c_0\varepsilon^2$, 
where $c_0$ is a positive absolute constant ???.
We may assume for the John ellipse of ${\tilde{K}}$
that $E_i(\widetilde{K})=B^2$.
Then there is a point $p$ of $\widetilde{K}$ of distance
at least  $1+c_0\varepsilon^2$ from $o$.
We deduce that $V(\widetilde{K}) \geq \pi+
c_0^{3/2} \left( 2 \sqrt{2}/3) + o(1) \right) \varepsilon ^3$, for
$\varepsilon \to 0$. Similarly, $({\tilde{K}})^*$ has a
supporting line $p^*$ at distance at most $1 - \left( c_0 + o(1) \right)
\varepsilon ^2$ from $o$, hence $V \left( {\tilde{K}} ^* \right)  \le \pi -
c_0^{3/2} \left( 4 \sqrt{2}/3) + o(1) \right) \varepsilon ^3$.
Therefore Theorem B yields 

\newpage

$$
\cases
4V({\tilde{K}})V \left( ({\tilde{K}})^* \right) = 
\left( V({\tilde{K}}) + V \left( ({\tilde{K}})^* \right) \right) ^2 -
\left( V({\tilde{K}}) - V \left( ({\tilde{K}})^* \right) \right) ^2 \\
\le 4 \pi ^2 - \left( \left( \sqrt{8} c_0^{3/2} + o(1) \right) \varepsilon ^3 
\right) ^2.
\endcases
$$
The $o$-symmetric case follows similarly.
$\blacksquare $
\enddemo


\definition{Remark}
A somewhat weaker result is given by 
K. M. Ball, K. J. B\"or\"oczky \cite{BB}, Theorem 2.2, case of ${\Bbb{R}}^2$ 
(observe that there $n \ge
3$ is written, but the arguments are valid for $n \ge 2$ as well). 
By a slight rewriting it
says the following. (Observe that in the end of the proof of the cited theorem
the exponent of $|\log \varepsilon |$ can in fact 
be doubled at the application of (3)
and (4) there, for the $o$-symmetric case.)
Under the conditions of Theorem D, except that $K \subset {\Bbb{R}}^n$,
we have
$$
V(K)V(K^*) \leq \kappa_n ^2 
(1-{\text{const}}_n\varepsilon^ {3n+3} | \log \varepsilon |^{-4}),
$$
and for the $o$-symmetric case 
$$
V(K)V(K^*) \leq \kappa _n ^2 
(1-{\text{const}}_n\varepsilon^ {(3n+3)/2} | \log \varepsilon |^{-4}).
$$
For $n=2$ these give 
$$
V(K)V(K^*) \leq \pi ^2 (1-{\text{const}}
\varepsilon^ 9 | \log \varepsilon |^{-4}),
$$
and
$$
V(K)V(K^*) \leq \pi ^2 (1-{\text{const}}
\varepsilon^ {9/2} | \log \varepsilon |^{-4}).
$$
As stated in \cite{BB}, we cannot have a better upper estimate, even in the
$o$-symmetric case, than 
$\kappa _n ^2 (1-{\text{const}}_n\varepsilon^ {(n+1)/2})$. 

For the planar $o$-symmetric case we can let $K_0:= \{ (x,y) \in B^2 \mid |x|
\le 1 - \varepsilon $. Then the
John ellipse $E_i(K)$ of $K$ has, by Behrends' characterization, the equation
$\left( x / (1 - \varepsilon ) \right) ^2 + y^2 = 1$. Possibly the
Banach-Mazur distance is attained for $E_i(K)$ and that inflated copy $
\lambda E_i(K)$ of it
that passes through the four non-smooth points $( \pm(1- \varepsilon , \pm
\sqrt{2 \varepsilon - \varepsilon ^2}$ of $K$. If this were true, we would
have $\lambda ^2 = \left( (1 - \varepsilon ) / (1 - \varepsilon ) \right) ^2
+ ( \sqrt{2 \varepsilon - \varepsilon ^2} )^2 $, hence $\lambda -1 \sim
\varepsilon $. Also we have $V(K)V(K^*)= \pi ^2 - 
\left( \pi \cdot 4 \sqrt{2} +o(1) \right) \varepsilon ^{3/2}/3$.
Observe that the same Banach-Mazur distance is attained if we consider the
intersection $K'$ of $B^2$ with an $o$-symmetric square $Q$ of sidelengths $1 -
\varepsilon $, that follows from the fact that the John and L\"owner ellipses
of $K'$ are the incircle of $Q$ and $B^2$, respectively. But $\left( 
V(K')V(K')^* \right) = \pi ^2 - 
\left( \pi \cdot 8 \sqrt{2} +o(1) \right) \varepsilon ^{3/2}/3$, so $K'$ is
worse than $K$. Similarly, for the general case the intersection of $B^2$ with
a regular triangle of sides at a distance $1 - \varepsilon $ from $o$ is
probably worse than the body $K:= \{ (x,y) \in B^2 \mid x
\le 1 - \varepsilon $. Whether $K_0$ and $K$ could be conjectured as optimal, 
is not clear: possibly some truncation of $B^2$ with some other curve(s) can
be better.

\enddefinition

\head About the stability of the Santal\'o point\endhead

\newpage

Let us review some basic formulas about polar bodies, due to
Santal\'o and  Meyer-
Pajor.
Let  $M$ be
a convex body in $R^d$. The support function of
$M$ is
$$
h_M(u)=\max_{x\in M} \langle u,x \rangle \,.
$$
In particular, if $z \in {\text{int}}\,M$, then the support function of
$M-z$ is $h_M(u)- \langle u,z \rangle $, and
the radial function of $(M-z)^*$ at $u \in S^{d-1}$
is $(h_M(u)- \langle u,z \rangle )^{-1}$. It follows that (with $V( \cdot )$
denoting volume)
$$
V((M-z)^*)=d^{-1}\int_{S^{d-1}} (h_M(u)- \langle z,u \rangle )^{-d}\,du. 
$$
Here $V((M-z)^*)$ is a strictly convex analytic function of 
$z \in {\text{int}}\, M$, 
that tends to infinity as dist\,$(z, \partial M) \to 0$. In fact, its second
differential is a positive definite quadratic form. E.g., 
$$
\left( \frac{ \partial } {\partial x_1 } \right) V((M-z)^*) =
(d+1) \int_{S^{d-1}} u_1^2 (h_M(u)- \langle z,u \rangle )^{-d}\,du >0, 
$$
and the respective foirmula holds for thecsecind derivetive along any
direction.
Hence it has a unique
minimum at some $s(M) \in {\text{int}}\, M$, 
which is called the Santal\'o point $s(M)$
of $M$. Differentiation yields that
$$
\int_{S^{d-1}} \frac{u}{(h_M(u)- \langle s(M),u \rangle )^{d+1}}\,du=o\,. 
$$
Thus $\int_{(M-s(M))^*}y\,dy=o$, and hence $o$
is the center of mass of $(M-s(M))^*$. This implies that
$$
-(M-s(M)) ^* \subset d (M-s(M))^*\,,
$$
hence
$$
-(M-s(M)) \subset d(M-s(M))\,.
\tag 0
$$


There are known several statements about stability of the Santal\'o point, or
behaviour of $V(K-z)^*$ for $z$ close to the Santal\'o point $s(K)$ of $K$.
See e.g., Santal\'o \cite{San}, {\bf{2}}, pp. 156-157,
Kim-Reisner \cite{KR}, Propositions 1 and 2, and the first part of this
paper \cite{BMMR}, Lemma 11, first statement (about $c_1(K_0)$ --- the second
statement, about $c_2(K_0)$, will not be needed here).

%

Now we cite 
\cite{BMMR}, Lemma 11, first part. 
%
\cite{BMMR}, Lemma 11 becomes correct.)
This is a more explicit version of
Kim-Reisner \cite{KR}, Proposition 1, inasmuch the constants are explicitly
given. This will render it possible to give estimates with 
absolute constants (for fixed dimension).


\newpage

\proclaim{Lemma E'} (\cite{BMMR}, Lemma 11, first part) 
Let $d \ge 2$ be an integer, 
$K_0 \subset {\Bbb R}^d$ be a convex body, 
and let
$0 < \varepsilon _1 \le \varepsilon _1(K_0):= \min \, \{ 1/2,
2^{-2d-1}
\left( \kappa _{d-1}/(d \kappa _d ^2) \right) \cdot $
$|K_0|
/({\text{\rm{diam}}}\,K_0)^d \}$. 
Let $K \subset {\Bbb R}^d$ be a convex body, and let 
$$
(1- \varepsilon _1) K_0 + a 
\subset K \subset (1+ \varepsilon _1)  K_0 +b, {\text{ where }} a,b \in 
{\Bbb R}^d.
$$ 
Additionally, let \thetag{E2} from above hold.
Then 
$$
\| s(K)-s(K_0) \| \le c(K_0) \cdot \varepsilon _1\,,
\tag E2.5
$$
where 
$$
c(K_0):=
({\text{\rm{diam}}}\,K_0)^{(d+1)^2} 
|K_0|^{-d-2} \cdot d (d \kappa _d / \kappa _{d-1}) ^{d+2} \,.
\tag E3
$$
\endproclaim


\proclaim{Theorem E} Let $d \ge 2$ be an integer. Then there exists an 
$\varepsilon (d) > 0$ such that 
for all $\varepsilon  \in [0, \varepsilon (d))$
the following holds.
Let $K_0,K \subset {\Bbb R}^d$ be convex bodies, and let 
$$
(1- \varepsilon ) K_0 + a 
\subset K \subset (1+ \varepsilon )  K_0 +b,
$$
where $a,b \in {\Bbb R}^d$
and 
$$
[ \left( (1- \varepsilon ) K_0 + a \right) + \left( (1+ \varepsilon )  K_0 +b
\right) ] / 2 = K_0.
$$
Then 
we have $s(K) \in {\text{{\rm{int}}}}\,K_0$, and
$$
\cases
V((K_0-s(K))^*) - V((K_0-s(K_0))^*) \le \\
2^{2d^2+4d+1}d^{2d^2+6d+9}\kappa _{d-1} ^{-2d-4} \kappa _d (d+1)^{d+2}
\cdot \varepsilon ^2.
\endcases
$$
\endproclaim


\demo{Proof}
We proceed on the lines of Santal\'o \cite{San} and Kim-Reisner \cite{KR},
Proposition 2 and \cite{BMMR}, Lemma 11, second part. 

Observe that the statement of Theorem E is affine invariant. Therefore we may
suppose, by John's theorem, that $B^d \subset K_0 \subset d B^d$. Then
\thetag{E3} remains valid 
if in the expression of $c(K_0)$ in Lemma E' we replace
$|K_0|$ with its lower bound $\kappa _d$ and diam\,$K_0$ by its upper bound
$2d$, obtaining
$$
\| s(K)-s(K_0) \| \le c \cdot \varepsilon \,,
\tag E4
$$
where
$$
c=2^{(d+1)^2} d^{d^2+3d+4} \kappa _{d-1} ^{-d-2}.
\tag E5
$$
 
The minimum of $f(z):=V(\left( K_0 - x ) ^* \right) $ 
for $z \in {\text{int}}\,K_0$ is attained for $z=s(K_0)$,
the Santal\'o point of $K_0$. The function $f(z)$ is analytic, hence has a power
series 

\newpage

expansion at $s(K_0)$, convergent
in any open ball of centre $s(K_0)$ and contained in int\,$K_0$.
Let us suppose that $s(K)-s(K_0)$ is non-zero (else we have nothing to prove),
and it points, say, to the direction of the positive $x_d$-axis.
The first derivatives of $f(z)$ vanish at $z=s(K_0)$, hence 
$$
\cases
V(\left( K_0 - s(K) ) ^* \right) - V(\left( K_0 - s(K_0) ) ^* \right) = \\
f \left( s(K) \right) - f \left( s(K_0) \right) = (1/2) \| s(K)-s(K_0) \| ^2
\cdot (\partial / \partial x_d)^2 f(x) = \\
(1/2) \| s(K)-s(K_0) \| ^2 \cdot (d+1) \int _{S^{d-1}} u_d^2
\left( h_{K_0} - \langle u,x \rangle \right) ^{-d-2} du ,
\endcases
\tag E6
$$
where $x$ lies in the relative interior of the segment $[s(K_0) ,s(K)]$, and
hence $\| x - s(K_0) \| < \| s(K) - s(K_0) \| $.
We estimate $(\partial / \partial x_d)^2 f(x)$ in \thetag{E6}.
We estimate $u_d^2$ from above by $1$. Then still we have to estimate 
$h_{K_0} - \langle u,x \rangle $ from below by some positive number.
We have
$$
h_{K_0}(u) - \langle u,x \rangle = \left( h_{K_0}(u) - \langle u, s(K_0) 
\rangle \right) + 
\left( \langle u, s(K_0) \rangle -\langle u, x \rangle \right) .
\tag E7
$$ 
So we have to estimate from below $\left( h_{K_0}(u) - \langle u, s(K_0) 
\rangle \right) $ and to estimate from above 
$ | \langle u, s(K_0) \rangle -\langle u, x \rangle | $.

For the first estimate we recall \thetag{0}: 
$$
K_0 - s(K_0) \subset -d \left( K_0 - s(K_0) \right) ,
$$
or, equivalently,
$$
h_{K_0 - s(K_0)} \le d h_{s(K_0) - K_0} .
\tag E8
$$
This last inequality 
can be interpreted geometrically as follows. Take any supporting parallel
strip of
$K_0$. Then for the (positive)
distances of $s(K_0)$ from the two boundary hyperplanes
of the parallel strip we have that their quotient is at most $d$, or,
otherwise said, the distance of $s(K_0)$ to any of these boundary hyperplanes
cannot be less than $1/(d+1)$ times the width of the parallel strip.
We apply this for $K_0$. Thus the distance of $s(K_0)$ to 
any supporting hyperplane of $K_0$ cannot be less than $1/(d+1)$ 
times the minimal
width of $K_0$, which width is by $K_0 \supset B^d$ at least $2$. 
Returning to the original formulation, we have that 
$$
h_{K_0} - \langle s(K_0), u \rangle \ge 2/(d+1) .
\tag E9
$$

For the second estimate we have, by $\| u \| = 1$, 
$$
\cases
| \langle u, x \rangle -\langle u, s(K_0) \rangle | =  \\
|\langle u, x - s(K_0) \rangle | \le \| x - s(K_0) \| < \| s(K) - s(K_0) \| 
\le c_1 \varepsilon < c_1 \varepsilon (d).
\endcases
\tag E10
$$

\newpage

Suppose that 
$$
c_1 \varepsilon (d) \le 1/(d+1).
\tag E11
$$ 
Then, using \thetag{E7}, \thetag{E9} and \thetag{E10}, we get 
$$
\cases
h_{K_0}(u) - \langle u,x \rangle = \left( h_{K_0}(u) - \langle u, s(K_0) 
\rangle \right) + \\
\left( \langle u, s(K_0) \rangle -\langle u, x \rangle \right) \ge 2/(d+1) -
1/(d+1) = 1/(d+1) .
\endcases
\tag E12
$$
This implies, by \thetag{E6}, \thetag{E12} and \thetag{E4}, \thetag{E5}, that
$$
\cases
V[\left( K - s(K) \right) ^* ] - V[\left( K - s(K_0) \right) ^* = 
(1/2) \| s(K) - s(K_0) \| ^2 \cdot \\
\int _{S^{d-1}} u_d^2   
\left( h_{K_0}(u) - \langle u,x \rangle \right) ^{-d-2} du \le \\
(1/2) c^2 \varepsilon ^2 \cdot d \kappa _d [1/(d+1)] ^ {-d-2} ,
\endcases
\tag E13
$$
and this is equivalent to the inequality of the Theorem.
\enddemo

\head Optimality of regular $n$-gons\endhead

We introduce some notation that will be used
in the next section, as well.
For $x \neq y\in R^2$, we 
write  aff$\, \{ x, y \} $ for
the line passing through $x$ and $y$.
For compact sets $X_1,\ldots,X_k$ in $R^2$, let
$[X_1,\ldots,X_k]$ denote the convex hull of their union.
For any point $p \neq o$, let $p^*$ be the polar line
with equation $\langle x, p \rangle =1$. For 
$0 \not\in l=p^*$, let $p=l^*$.


\proclaim{Theorem F} Among convex polygons $K$ of at
most $n$ vertices whose Santal\'o point is $o$, the ones
maximizing the area product $V(K)V(K^*)$ are the non-singular
linear images of regular
$n$-gons with centre at $o$.
\endproclaim

\demo{Proof} As all triangles are affine images of each other, 
we may assume $n \geq 4$. 
By Theorem ? we may assume even $n \ge 5$. Observe that 
the area product of a regular
$n$-gon, i.e., $n^2 \sin ^2 ( \pi /n )$ is strictly increasing with
$n$. Therefore,
using induction w.r.t. $n$, with base of induction $n=4$, we may restrict our
attention to convex polygons with exactly $n$ sides. 
By the intersection point
of two distinct parallel lines, we mean their common point of infinity (in the
projective plane). 

Let  $K$ maximize
the area product among convex polygons of at
most $n$ vertices, whose Santal\'o point is $o$. 
We prove two basic properties of $K$:

\medskip

\noindent(i) If $x_1,x_2,x_3$ are consecutive vertices of $K$,
then $o$, $\frac12(x_1+x_3)$ and $x_2$ are collinear. 

\medskip

\noindent(ii) If $x_1,x_2,x_3,x_4$ are consecutive vertices of $K$,
and $m$ is the intersection point of aff$\, \{ x_1, x_2 \} $
and aff$\, \{ x_3, x_4 \} $,
then $o$, $(x_2+x_3)/2$ and $m$ are collinear. Here $m$ is meant as a point of
the projective plane.

\medskip

To prove these claims, we show with a method going back to Mahler that 
if either (i) or (ii)
does not hold, then $K$ can be deformed into a convex polygon with the
same number of vertices but with a larger area product.
We write $\gamma _1,\gamma _2,\ldots$ 

\newpage

for
{\it{positive constants that depend on $K$, but not on its deformations}}.

First we suppose that (i) does not hold, and seek a contradiction.
For $i=1,3$, 
we write
$y_{1+i/2}$ for the intersection of $x_i^*$ and $x_2^*$,
and hence $y_{1.5}$ and $y_{2.5}$ are consecutive vertices of $K^*$. 
Let $l$ be the line
through $x_2$ parallel to aff$\, \{ x_1, x_3 \} $,
and hence $l^* \in [y_{1.5},y_{2.5}]$.
We consider two cases depending on the
position of $o$ with respect to the triangle $[x_1,x_2,x_3]$.

\medskip
 

\noindent{\bf Case 1. } $o \not\in {\text{int}}\,[x_1,x_2,x_3]$

We may assume that $x_1$ and $x_2$ lie on the
same side of aff$\,\{ o,(x_1+x_3)/2 \} $ if
$o \not\in [x_1,x_2,x_3]$ (and hence $o \not\in [x_1x_2x_3]$), 
or $\| x_1 \| < \| x_3 \| $
if $o \in [x_1,x_3]$. In particular,
$$
\| l^*-y_{1.5} \| > \| l^*-y_{2.5} \| \,.
$$

We move $x_2$
into a position $\tilde{x}_2$ where
$\tilde{x}_2-x_2=\varepsilon(x_3-x_1)$ for small $\varepsilon >0$,
where $\varepsilon $ is small enough to ensure that
$x_1$ and $\tilde{x}_2$ lie on the
same side of the line containing the other side of $K$ 
with endpoint $x_3$, and if $o \not\in [x_1,x_2,x_3]$,
then still $x_1$ and $\tilde{x}_2$ lie on the
same side of  aff$\, \{ o,(x_1+x_3)/2 \} $.
Therefore there exists 
a convex polygon $\widetilde{K}$ obtained from $K$ by removing
$x_2$ from the set of vertices, and adding $\tilde{x}_2$.
Clearly $V(\widetilde{K})=V(K)$.

For $i=1,3$, let
$\tilde{y}_{(i+2)/2}$ be the intersection of $x_i^*$ and $\tilde{x}_2^*$,
and hence $\tilde{y}_{1.5}$ and $\tilde{y}_{2.5}$ are
the two new vertices of $\widetilde{K}^*$ 
replacing the vertices $y_{1.5}$ and $y_{2.5}$ of $K^*$.
In addition, $\tilde{y}_{1.5} \not\in K^*$ and $\tilde{y}_{2.5}\in K^*$,
because $o \not\in {\text{int}}\,[x_1,x_2,x_3]$. Writing $\alpha :=\angle
x_2o{\tilde x}_2$ for
the angle of $x_2$ and $\tilde{x}_2$???, we have
$\gamma _1 \varepsilon < \alpha < \gamma_2 \varepsilon $.
Since
$$
V(\tilde{K}^*)-V(K^*)=V([l^*,y_{1.5},\tilde{y}_{1.5}])-
V([l^*,y_{2.5},\tilde{y}_{2.5}]),
$$
and $\alpha $ is the angle of both triangles $[l^*,y_{1.5},\tilde{y}_{1.5}]$
and $[l^*,y_3,\tilde{y}_3]$ at $l^*$???, we have
$$
V(\tilde{K}^*)-V(K^*) > \frac{(\| l^*-y_{1.5} \| ^2-\| l^*-y_{2.5} \| ^2)
\sin \alpha}{2}-\gamma _3 \alpha ^2 >
\gamma _4 \varepsilon \,.
$$
It follows by Theorem E that 
$V((\tilde{K}-s(\tilde{K}))^*) > V(\tilde{K}^*)-\gamma _5 \varepsilon ^2$.
Therefore if $\varepsilon >0$ is small enough, 
then $V \left( \left( {\tilde{K}}-s({\tilde{K}}) \right) ^* \right) 
> V(K^*) + \gamma _6 \varepsilon > V(K^*)$, 
which is a contradiction, and thus we have
proved (i) in Case 1.

\medskip

\noindent{\bf Case 2. } $o \in {\text{int}}\,[x_1,x_2,x_3]$

We may assume that $x_1$ and $x_2$ lie on the
same side of aff$\,\{ o,(x_1+x_3)/2 \} $, and hence
$$
\| l^*-y_1 \| < \| l^*-y_3 \|\,.
$$

We move $x_2$
into a position $\tilde{x}_2$ where
$\tilde{x}_2-x_2=\varepsilon (x_1-x_3)$ for small $\varepsilon >0$,
and $\varepsilon $ is small enough to ensure that
there exists 
a convex polygon $\widetilde{K}$ obtained from $K$ by removing
$x_2$ from the set of vertices, and adding $\tilde{x}_2$.
Clearly $V(\widetilde{K})=V(K)$.

For $i=1,3$, let again
$\tilde{y}_{(i+2)/2}$ be the intersection of $x_i^*$ and $\tilde{x}_2^*$,
and hence $\tilde{y}_{1.5}$ and $\tilde{y}_{2.5}$ are
the two new vertices of $\widetilde{K}^*$ 
replacing the vertices $y_{1.5}$ and $y_{2.5}$ of $K^*$.
In this case $\tilde{y}_{1.5} \in K^*$ and $\tilde{y}_{2.5} \not\in K^*$,
thus

\newpage

$$
V(\tilde{K}^*)-V(K^*)=V([l^*,y_{2.5},\tilde{y}_{2.5}])-
V([l^*,y_{1.5},\tilde{y}_{1.5}]).
$$
Now the argument can be finished as in Case 1, completing
the proof of (i).

\medskip

Next we suppose that (ii) does not hold.
We may assume that 
$$
\cases
[x_1,x_2] {\text{ and }} p=(x_2+x_3)/2  {\text{ lie on the same side}} \\
{\text{of the line connecting }} o {\text{ and }} m\,. 
\endcases
\tag *
$$
There are three cases. The point $m$ can be a finite point, lying on the other
side of aff\,$\{ x_2,x_3 \} $ as $K$, or can be a point at infinity, or can be
a finite point lying on the same side of aff\,$\{ x_2,x_3 \} $ as $K$.
Let $l_{1.5}={\text{aff}}\, \{ x_1, x_2 \} $, $l_{3.5}=
{\text{aff}}\,\{ x_3, x_4 \} $ and $l={\text{aff}}\,\{ x_2, x_3 \} $.
We observe that if $m\in {\Bbb{R}}^2$, then 
$\langle m,l_{1.5}^* \rangle = \langle m,l_{3.5}^* \rangle =1$, and also 
$\langle 0,l_{1.5}^* \rangle = \langle 0,l_{3.5}^* \rangle =0$.
Therefore the assumption on $p$ yields that 
$$
\langle p, l_{1.5}^* \rangle > \langle p,l_{3.5}^* \rangle
\tag **
$$ 
holds independently whether $m\in {\Bbb{R}}^2$, or $m$ is a point at infinity.

We rotate $l$ about $p$ into a new position $\tilde{l}$ through a small angle 
$\varepsilon  \in (0,\pi/2)$,
in a way such that $\tilde{l}$ intersects $[x_3,x_4]\setminus \{ x_4 \} $ in
a point $\tilde{x}_3$, but intersects $l_{1.5}$ in a point ${\tilde x}_2$
lying outside of $K$. In particular,
$$
\| \tilde{l}^*-l^* \| > \gamma _6 \varepsilon \,.
$$
Let $\widetilde{K}$ be the convex polygon obtained from $K$ by removing
$x_2$  and $x_3$ from the set of vertices, 
and adding $\tilde{x}_2$ and  $\tilde{x}_3$.
Then
$$
V(\tilde{K})=V(K)+V([p,x_2,\tilde{x}_2])-
V([p,x_3,\tilde{x}_3]) > V(K)-\gamma _7 \varepsilon ^2 \,.
$$
Now $\widetilde{K}^*$ is obtained from $K^*$ by removing
$l^*$ from the set of vertices, and adding $\tilde{l}^*$.
Here $l^*,\tilde{l}^*\in p^*$, 
in all the three cases for $m$, i.e., if it is a point at infinity, or it is a
finite point, lying on either side of aff\,$\{ x_2,x_3 \} = l$. We have three
different figures, but the following calculations are valid for each of these
three figures. Supposing that the vertices of $K$ follow each other
in the positive sense, then
assumption {\thetag{*}}
on $p$ yields $\langle ( l_{3.5}^*-l_{1.5}^*) \times ({\tilde l}^*-l^*) ,
(0,0,1) \rangle >0$, hence by {\thetag{**}}
actually
$$
\langle ( l_{3.5}^*-l_{1.5}^*) \times ({\tilde l}^*-l^*) , (0,0,1)
\rangle > \gamma _8 \varepsilon 
$$
(observe that the angle of $m^*$ and $p^*$ only depends on $K$, but not on
$\varepsilon $).
Therefore

\newpage

$$
V(\tilde{K}^*)-V(K^*)=\langle 
( l_{3.5}^*-l_{1.5}^*) \times ({\tilde l}^*-l^*) , (0,0,1) \rangle /2
> \gamma _9 \varepsilon \,.
$$
Furthermore, 
$V((\tilde{K}-s(\tilde{K}))^*) > V(\tilde{K}^*)-\gamma _{10} \varepsilon ^2$
by Theorem E.
We conclude that
$$
\cases
V(\tilde{K})V((\tilde{K}-s(\tilde{K}))^*)
> \left( V(K)-\gamma _7 \varepsilon ^2 \right) \left( V(K^*) - \gamma _{10}
\varepsilon ^2 \right) \\
> \left( V(K)-\gamma _7 \varepsilon ^2 \right) 
\left( V(K^*)+\gamma _9 \varepsilon - \gamma _{10} \varepsilon ^2 \right) \\
> V(K)V(K^*) + \gamma _{11} \varepsilon 
> V(K)V(K^*)\,,
\endcases
$$
provided $\varepsilon \in (0, \pi /2)$ 
is small enough. This contradiction proves (ii).

Now we prove Theorem F based on (i) and (ii).
Let $x_1, \ldots , x_k$, with $3 \leq k \leq n$ 
be the vertices of $K$ in this order. By the beginning of the proof of this
theorem we have $k=n \ge 5$. By $n \ge 5$ we have that the average value of 
$\angle x_{i-1}ox_i + \angle x_iox_{i+1}$ is $4 \pi /n < \pi$ (indices meant
cyclically). Let, e.g., $\angle x_1ox_2 + \angle x_2ox_3 < \pi $. (Observe
that this property is invariant under linear maps.)
Applying a linear transformation, we may assume that $x_1,x_2,x_3 \in S^1$
and hence $x_2$ lies on the
perpendicular bisector of $[x_1,x_3]$ by (i). 
Applying another linear transformation,
we may also assume that $x_1,x_2,x_3 \in S^1$  --- actually, they lie in an open
half-circle.

Now (ii) yields that aff$\, \{ x_1, x_2 \} $
and aff$\, \{ x_3, x_4 \} $ are symmetric with respect to
the perpendicular bisector of $[x_2,x_3]$, and $m$ is a finite point separated
from $K$ by $l$. It follows
that aff$\, \{ x_2, x_3 \} $
and aff$\, \{ x_3, x_4 \} $ are symmetric with respect to
aff$\,\{ o,x_3 \} $. Therefore (i) yields that
$x_4 \in S^1$, with $\| x_4-x_3 \| = \| x_3-x_2 \| = \| x_2-x_1 \| $.
Continuing like this, we conclude that
$K$ is a regular $k$-gon
inscribed into $B^2$. From the first paragraph of the proof we have $k=n$.
$\blacksquare $
\enddemo



\Refs

\widestnumber\key{WWW}



\ref
\key BB
\by Ball, K. M., B\"or\"oczky, K. J.
\paper Stability of some versions of the Pr\'ekopa-Leindler inequality
\jour Monatshefte Math.
\vol 163
\yr 2011
\pages 1-14
\MR {\bf{2012b:}}{\rm{52015}} 
\endref 

\ref
\key B 
\by Behrends
\paper \"uber die kleinste umbeschriebene und die gr\"osste einbeschriebene
Ellipse eines konvexen Berichs
\jour Math. Ann.
\vol 113
\yr 1937
\pages 379-411
Zbl. {\bf{18.}}17502
\endref 

\ref 
\key J 
\by John
\paper Extremum problems with inequalities as subsidiary conditions
\jour 
Studies and Essays Presented to R. Courant on his 60th Birthday, Interscience,
New York, N. Y., 1948
\pages 187-204
\endref   




\endRefs
\enddocument